% PlainTeX
\def\vers{Sept.~30, 2008, v.1}
\magnification=1200
\hsize=6.5truein
\vsize=8.9truein
\font\bigfont=cmr10 at 14pt
\font\mfont=cmr9
\font\sfont=cmr8
\font\mbfont=cmbx9
\font\sifont=cmti8
\def\scirc{\,\raise.2ex\hbox{${\scriptstyle\circ}$}\,}

\def\a{\alpha}
\def\C{{\bf C}}
\def\D{\Delta}
\def\d{\delta}

\def\g{\gamma}
\def\G{\Gamma}
\def\GG{\widehat{\Gamma}}
\def\H{{\bf H}}

\def\h{\hbox}

\def\L{{\cal L}}
\def\LL{\widehat{\cal L}}

\def\n{\nu}

\def\O{{\cal O}}
\def\Q{{\bf Q}}
\def\q{\quad}
\def\qq{\qquad}
\def\R{{\bf R}}
\def\S{\hskip1pt\overline{\!S}\hskip1pt}

\def\s{\sigma}
\def\V{{\cal V}}
\def\VV{\widehat{\cal V}}
\def\Z{{\bf Z}}

\def\Ext{\hbox{\rm Ext}}
\def\Hom{\hbox{\rm Hom}}
\def\Ker{\hbox{\rm Ker}}
\def\Gr{\hbox{\rm Gr}}

\def\MHS{{\rm MHS}}

\def\Im{{\rm Im}}
\def\NF{{\rm NF}}
\def\ad{{\rm ad}}
\def\1{{\hskip1pt}}

\def\into{\hookrightarrow}
\hbox{}
\vskip 1cm

\centerline{\bigfont A Generalization of the N\'eron Models}

\medskip
\centerline{\bigfont of Green, Griffiths and Kerr}

\bigskip
\centerline{Patrick Brosnan, Gregory Pearlstein, and Morihiko Saito}

\bigskip\medskip
{\narrower\noindent
{\mbfont Abstract.} {\mfont
We generalize a construction of the N\'eron model for a family
of intermediate Jacobians due to Green, Griffiths and Kerr
by using the theory of mixed Hodge modules.
It is a topological group defined over any partial compactification
of the base space, and it `graphs' admissible normal functions.
Moreover, there is a stratification of the partial compactification
such that the restriction over each stratum is a complex Lie group
over the stratum.}
\par}

\bigskip\bigskip
\centerline{\bf Introduction}
\footnote{}{{\sifont Date\1}{\sfont:\ \vers}}

\bigskip\noindent
Let $\H$ be a polarizable variation of $\Z$-Hodge structure of
weight $w<0$ on a complex manifold $S$.
Let $J_S(\H)$ be the family of intermediate Jacobians over $S$.
Its fiber at $s\in S$ is
$$J_S(\H)_s=\L(s)/(F^0\L(s)+L_{\Z,s}),$$
where $(\L,F)$ is the underlying filtered locally free sheaf of
$\H$, and $L_{\Z}$ is the underlying local system of $\H$.
Let $j:S\to\S$ be a partial compactification as a complex analytic
space.
Let $\{S_{\a}\}$ be a Whitney stratification of $\S$ such that
$S$ is one of the strata.

\medskip\noindent
{\bf Theorem.} {\it There is a N\'eron model $J_{\S}(\H)$ over
$\S$ which extends $J_S(\H)$ and such that any admissible normal
function on $S$ uniquely extends to a section of $J_{\S}(\H)$.
Moreover, there is a short exact sequence of topological groups
over $\S$
$$0\to J_{\S}(\H)^0\to J_{\S}(\H)\to G\to 0,$$
such that the restriction of $J_{\S}(\H)^0$ over $S_{\a}$ is
a complex Lie group over $S_{\a}$ with connected fibers and
$G_s$ is discrete and is a subgroup of $(R^1j_*L_{\Z})_s$ for
any $s\in\S$.}

\medskip
In the curve case this was proved first by Clemens [3] in some
special cases, and then by [10] in the general curve case where it
was defined by enlarging the Zucker extension [13].
However, it was pointed out by Green, Griffiths and Kerr [7]
that a subspace of it is sufficient and is more natural (still in
the curve case).
In this paper we show that their construction can be extended
naturally to the general case using the theory of mixed Hodge
modules.

In Section 1 we recall some basics about Zucker extensions and
admissible normal functions.
In Section 2 we prove Theorem by constructing the N\'eron models.

\bigskip\bigskip
\centerline{\bf 1.\ Zucker Extensions and Admissible Normal
Functions}

\bigskip\noindent
{\bf 1.1.~Zucker extensions.}
Let $\H=((\L,F),L_{\Z})$ be a polarizable variation of $\Z$-Hodge
structure of weight $w<0$ on a complex manifold $S$.
Here $(\L,F)$ is the underlying filtered locally free sheaf and
$L_{\Z}$ is the underlying $\Z$-local system.
We assume that $L_{\Z}$ is torsion-free in this paper.
Let $\V$ be the vector bundle on $S$ corresponding to the locally
free sheaf $\L/F^0\L$, and let $\G\subset\V$ denote the subgroup
over $S$ corresponding to the subsheaf $L_{\Z}\subset\L/F^0\L$
where the last injectivity follows from the negativity of the
weight. Set
$$J_S(\H)=\V/\G.$$
Here the quotient is set-theoretically taken fiberwise.
This has a structure of a complex Lie group over $S$.

Let $j:S\to\S$ be a smooth partial compactification of $S$ such that
$D:=\S\setminus S$ is a divisor with normal crossings.
Assume the local monodromies of $L_{\Z}$ are unipotent.
Let $\LL$ be the Deligne extension of $\L$, see [4].
By Schmid [11], the Hodge filtration $F$ on $\L$ is uniquely extended
to a filtration $F$ on $\LL$ such that $\Gr_F^p\LL$ are locally free
(i.e. $F^p\LL=\LL\cap j_*F^p\L$).
Let $\VV$ denote the vector bundle on $\S$ corresponding to the
locally free sheaf $\LL/F^0\LL$.
Let $\GG$ be the subgroup of $\VV$ over $\S$ corresponding to
$j_*L_{\Z}\subset\LL/F^0\LL$.
Then the Zucker extension is defined by
$$J_{\S}^Z(\H)=\VV/\GG.$$

\medskip\noindent
{\bf 1.2.~Restriction to the diagonal curve.}
Let $C$ be the diagonal curve defined by $z_i=z_j$ for $i\ne j$,
where the $z_i$ are local coordinates such that
$S=\bigcap\{z_i\ne 0\}$.
Since the local monodromies are unipotent, the restriction of
the Deligne extension to $C$ is again the Deligne extension.
This follows from the fact that the restriction of
$\exp(-\sum_i(\log z_i)N_i)u$ to the diagonal curve is
$\exp(-(\log z)\sum_iN_i)u$ where $u$ is a multivalued horizontal
section.
(Here $N_i=(2\pi i)^{-1}\log T_i$ with $T_i$ the local monodromies.)

As a corollary, we see that the restriction of the Zucker extension
to the diagonal curve is again the Zucker extension.

\medskip\noindent
{\bf 1.3.~Admissible normal functions.}
A normal function is a holomorphic section $\n$ of $J_S(\H)$ satisfying
the Griffiths transversality.
By Carlson [2] it corresponds to an extension class of $\Z_S$ by
$\H$ giving a short exact sequence
$$0\to\H\to\H'\to\Z_S\to 0.\leqno(1.3.1)$$
A normal function $\n$ is called admissible with respect to a partial
compactification $\S$ of $S$ if $\H'$ is an admissible variation of
mixed Hodge structure ([8], [12]) with respect to $\S$, see [10].
The group of such normal function will be denoted by
$\NF(S)^{\ad}_{\S}$.

Assume $\S$ is as in (1.1), i.e. $D$ is a divisor with norma
crossings, and moreover the local monodromies are unipotent.
By [8] (and [12] in the 1-dimensional case) the conditions
for an admissible normal function are given as follows:
\smallskip
(i) $\Gr_F^p\Gr_k^W\LL'$ are locally free.

\smallskip
(ii) The relative monodromy filtration exists for any local
monodromy.

\smallskip\noindent
For $s\in\S$, we have the cohomological invariant of $\n$ at $s$
$$\g_s(\n)\in H^1(U_s,L_{\Z}|_{U_s})=
\Ext^1(\Z_{U_s},L_{\Z}|_{U_s}),$$
where $U_s$ is the intersection of $S$ with a sufficiently small
ball with center $s$ in an ambient space.
This invariant is defined by passing to the underlying short exact
sequence of $\Z$-local systems, and restricting it over $U_s$,
see also [1].

In the curve case it is known that the normal function
extends to a section of the Zucker extension if the
cohomological invariant $\g_s(\n)$ vanishes, see e.g. [10],
Prop.~2.3 and also [6] for the geometric case.
Indeed, if $\g_s(\n)=0$, then $\n$ is given by the class of
$$\s_F(1)-\s_{\Z}(1)\in\G(\S,\LL),\leqno(1.3.2)$$
where $\s_{\Z},\s_F$ are splittings of the
underlying short exact sequence of locally free sheaves of (1.3.1)
$$0\to\LL\to\LL'\to\O_{\S}\to 0,$$
such that $\s_{\Z}$ is defined over $\Z$ and $\s_F$ is compatible
with $F$.
This argument can be extended to the normal crossing case of
higher dimension.
(In the curve case its converse is also true, see [10], Prop.~2.4.)

\medskip\noindent
{\bf 1.4.~Relation with the limit mixed Hodge structure.}
In the curve case, assuming the local monodromy is unipotent,
the limit mixed Hodge structure $\psi_t\H$ at $0\in\S\setminus S$
is given by
$$\eqalign{H&=((H_{\C};F,W),(H_{\Q},W),H_{\Z}),\cr
\h{with}\qq H_{\C}&=\LL(0)\,(:=\LL/m_0\LL),\q
H_A=\psi_tL_A\,\,\,(A=\Z,\Q),}$$
where $m_0\subset\O_{\S,0}$ is the maximal ideal and $t$ is a
local coordinate of $\S$.
Note that the nearby cycle functor $\psi_tL_A$ can be defined
in this case by
$$\psi_tL_A=\Gamma(\widetilde{\D}^*,\rho^*L_A),$$
where $\D$ is a sufficiently small disk around $0$ and
$\rho:\widetilde{\D}^*\to\D^*$ is a universal covering.

We have a commutative diagram of MHS
$$\matrix{0&\to&\psi_t\H&\to&\psi_t\H'&\to&\Z&\to&0\cr
&&\,\,\,\downarrow{\scriptstyle\!N}&\raise12pt\hbox{ }
\raise-8pt\hbox{ }&\,\,\,\downarrow{\scriptstyle\!N}&&
\,\,\,\downarrow{\scriptstyle\!0}\cr
0&\to&\psi_t\H(-1)&\to&\psi_t\H'(-1)&\to&\Z(-1)&\to&0\cr}$$
The horizontal exact sequence defines
$$\n'_0\in J(\psi_t\H).$$
In the case $\g_0(\n)=0$, let $\n''_0$ denote the value at 0 of the
normal function given by (1.3.2). Then we have by the definition of
the limit mixed Hodge structure (and using [2])
$$\n'_0=\n''_0\in J(\psi_t\H).$$
Set $H_0=\Ker(N:\psi_t\H\to\psi_t\H(-1))$.
(This is compatible with the definition of $H_s$ in (2.1) below.)
By the snake lemma applied to the above diagram, there is
$$\n_0\in J(H_0),$$
whose injective image in $J(\psi_t\H)$ coincides with $\n'_0$,
since the connecting morphism gives $\g_0(\n)$.
Note that this implies another proof of a theorem of Green,
Griffiths and Kerr [7], Thm.~(II.A.9)(i).
This $\n_0$ coincides with the one constructed in (2.1)
below in a more general situation, since the vertical arrow of the
above diagram calculates the functor $i_0^*\R j_*$ in this case.

\bigskip\bigskip
\centerline{\bf 2.\ N\'eron Models}

\bigskip\noindent
{\bf 2.1.~Identity components.}
Let $\n\in\NF(S)^{\ad}_{\S}$, i.e. be an admissible normal function
on $S$ with respect to a partial compactification $j:S\to\S$.
We have the corresponding short exact sequence
$$0\to\H\to\H'\to\Z_S\to 0.$$
For $s\in\S$, let $i_s:\{s\}\to\S$ denote the inclusion, and set
$$H_s:=H^0i_s^*\R j_*\H,\q J(H_s):=\Ext_{\MHS}^1(\Z,H_s),$$
where MHS denotes the category of mixed $\Z$-Hodge structures [5].
Note that the functors $H^ki_s^*\R j_*\H$ are defined in MHS
(with $\Z$-coefficients).
In fact, they are defined with $\Q$-coefficients in the algebraic
case [9],
and a similar argument works in the analytic case assuming
$\S\setminus S$ is an intersection of divisors (shrinking
$\S$ if necessary).
Moreover, forgetting the mixed Hodge structure, they are naturally
defined with $\Z$-coefficients so that
$$H^ki_s^*\R j_*L_{\Z}=
\mathop{\smash{\mathop{\h{\rm lim}}\limits_{\raise.4ex
\h{${\scriptstyle\longrightarrow}$}}}}
\limits_{\smash{\raise-0.7ex\h{$\scriptstyle U$}}}
H^k(U\cap S,L_{\Z}),\leqno(2.1.1)
$$
where $U$ runs over open neighborhoods of $s$ in $\S$.

We have an exact sequence of mixed $\Z$-Hodge structures
$$0\to H_s\to H^0i_s^*\R j_*\H'\to\Z\buildrel{\d_s}\over\to
H^1i_s^*\R j_*\H,$$
and the image of $1\in\Z$ by $\d_s$ is the cohomological invariant
$\g_s(\n)$ of the admissible normal function $\n$ at $s$.
So the exact sequence implies that, if $\g_s(\n)=0$, then $\n$
defines
$$\n_s\in J(H_s).\leqno(2.1.2)$$
We define set-theoretically the identity component by
$$J_{\S}(\H)^0:=\coprod_{s\in\S}J(H_s).$$
To define a topological structure on it, take a resolution of
singularities $\pi:\S'\to\S$ such that the pull-back of
$\S\setminus S$ is a divisor with normal crossings.
For $s'\in\S'$, we have the canonical injection of mixed Hodge
structures
$$\pi^*:H_{\pi(s')}\into H_{s'}.\leqno(2.1.3)$$
By (2.1.1) this is induced by the restriction morphisms
$$H^0(U\cap S,L_{\Z}) \into H^0(U'\cap S,L_{\Z}),$$
where $U,U'$ are open neighborhoods of $s,s'$ in $\S,\S'$ such
that $\pi(U')\subset U$.
To show that this respects the mixed Hodge structures,
set $Z=\pi^{-1}(s)$, $\pi_Z=\pi|_Z:Z\to\{s\}$, and let
$i_Z:Z\to \S'$, $i'_{s'}:\{s'\}\to Z$, $i_{s'}:\{s'\}\to\S'$,
$j':S\to\S'$ denote the inclusion morphisms.
Then (2.1.3) coincides with the composition of canonical morphisms
$$i_{s}^*\R j_*\H=i_{s}^*\R\pi_*\R j'_*\H=\R(\pi_Z)_*i_Z^*\R j'_*\H
\to i_{s'}^{\prime *}i_Z^*\R j'_*\H=i_{s'}^*\R j'_*\H.$$

The above description of (2.1.3) implies that the cokernel of
(2.1.3) is torsion-free since $H^0(U\cap S,L_{\Z})$ is identified
with the monodromy invariant subspace and $L_{\Z}$ is torsion-free.
Combining this with the negativity of the weight $w$, we see that
(2.1.3) induces an injection
$$J(H_{\pi(s')})\into J(H_{s'}).$$
Thus there is a subspace
$$J_{\S',\S}(\H)^0:=\coprod_{s'\in\S'}J(H_{\pi(s')})\subset
J_{\S'}(\H)^0=\coprod_{s'\in\S'}J(H_{s'}),$$
together with the surjection
$$J_{\S',\S}(\H)^0\to J_{\S}(\H)^0.$$
So the problem of constructing the N\'eron model is reduced to the
normal crossing case using the the quotient topology.
It is further reduced to the normal crossing case with
unipotent monodromy by taking locally a finite covering.
Then $J_{\S}(\H)^0$ is a subset of the Zucker extension, and we
have the induced topology on it.
So the topology is defined.
This construction is independent of the choice of the resolution.

Now let $\{S_{\a}\}$ be a Whitney stratification of $\S$ such that
the $R^kj_*L_{\Z}$ are locally constant on each stratum.
Let $i_{\a}:S_{\a}\to\S$ denote the inclusion.
Then $H^ki_{\a}\R j_*\H$ are variation of mixed Hodge structures.
So we have a structure of a complex Lie group on
$$J_{\S}(\H)^0|_{S_{\a}}=\coprod_{s\in S_{\a}}J(H_s).$$
Its underlying topology coincides with the induced
topology of the topology on $J_{\S}(\H)^0$ constructed above.
Indeed, we can show that an open subset $V$ of
$J_{\S}(\H)^0|_{S_{\a}}$ as a complex Lie group over $S_{\a}$ is an
open subset in the induced topology of the topology constructed
above as follows.
The assertion is reduced to the case $\S$ is smooth,
$\S\setminus S$ is a divisor with normal crossings, and the local
monodromies are unipotent.
Then there is an open subset $V'$ of $\GG$ such that its
restriction over $S_{\a}$ is the inverse image of $V$ to
$\VV|_{S_{\a}}$.
Taking the union of $V'$ with translates of $V'|_U$ by
local sections of $\GG$ defined over open subsets $U$ of $\S$,
we may assume that $V'$ is stable by the action of local sections
of $\GG$.
So the assertion follows.
(The other direction is easy.)

\medskip\noindent
{\bf 2.2.~Theorem.} {\it
Let $\n$ be an admissible normal function with respect to a partial
compactification $j:S\into\S$.
Assume $\g_s(\n)=0$ for any $s\in\S\setminus S$.
Then the $\n_s$ in $(2.1.2)$ define a continuous section of
$J_{\S}(\H)^0$, which is holomorphic over each stratum $S_{\a}$.}

\medskip\noindent
{\it Proof.}
By the definition of the topology, the assertion is reduced to the
normal crossing case with unipotent local monodromies.
Then it is sufficient to show that $\n_s$ coincides with the
section of the Zucker extension over $\S$ defined by (1.3.2)
using two sections $\sigma_F$ and $\sigma_{\Z}$.
So the assertion is reduced to the curve case by restricting to the
diagonal curves since the restriction of the Zucker extension to the
diagonal curve is again the Zucker extension, see (1.2).
In the curve case, we may assume $S=\D^*$ and $\S=\D$.
Then the assertion follows from (1.4).

\medskip\noindent
{\bf 2.3.~N\'eron models.}
Let $\n\in\NF(U\cap S)^{\ad}_U$ where $U$ is an open subset of $\S$.
It defines the $\n$-component $J_U(\H_{U\cap S})^{\n}$ of the
N\'eron model over $U$, which `graphs' $\n$.
More precisely, it has a canonical isomorphism
$$J_U(\H_{U\cap S})^{\n}=J_U(\H_{U\cap S})^0,\leqno(2.3.1)$$
such that $\n$ corresponds to the zero section of the right-hand
side.

For $\n,\n'\in\NF(U\cap S)^{\ad}_U$ such that
$\g_s(\n)=\g_s(\n')$ for any $s\in U\setminus S$, we have a
canonical isomorphism
$$J_U(\H_{U\cap S})^{\n}=J_U(\H_{U\cap S})^{\n'},\leqno(2.3.2)$$
which corresponds by the isomorphism (2.3.1) to an automorphism of
$J_U(\H_{U\cap S})^0$ defined by $\n-\n'$.
So the N\'eron model $J_{\S}(\H)$ is defined by gluing the
$J_U(\H_{U\cap S})^{\n}$ for admissible normal functions $\n$
defined over $U\cap S$ where $U$ are open subsets of $\S$.
This is done by induction on strata of a Whitney stratification of
$\S$ such that $R^1j_*L_{\Z}$ is a local system on each stratum
$S_{\a}$.
More precisely, we extend it over the open subsets
$\coprod_{\dim S_{\a}\ge r}S_{\a}$ by decreasing induction on $r$,
using the fact that $\g_s(\n)$ is locally constant on each stratum
$S_{\a}$.
Then the restriction of $J_{\S}(\H)$ over each stratum is a complex
Lie group over the stratum.
Here we use the fact that if $\g_s(\n)=\g_s(\n')$ then
$\g_{s'}(\n)=\g_{s'}(\n')$ for $s'\in\S\setminus S$
sufficiently near $s$.
For $s\in\S$, set
$$G_s=\bigcup_{U\ni s}\Im(\NF(U\cap S)^{\ad}_{U}\to
(R^1j_*L_{\Z})_s),\leqno(2.3.3)$$
where $U$ runs over open neighborhoods of $s$ in $\S$.
Then the topological group $G=\coprod_sG_s$ is identified with an
open subspace of $R^1j_*L_{\Z}$ viewed as an etale space associated
with a sheaf over $\S$.
As a conclusion, we get

\medskip\noindent
{\bf 2.4.~Theorem.} {\it The N\'eron model $J_{\S}(\H)$ graphs any
admissible normal functions $\n$.
More precisely, $\n$ defines a continuous section of
$J_{\S}(\H)$, which is holomorphic over each stratum $S_{\a}$ of
a Whitney stratification.
Furthermore, there is a short exact sequence of commutative
topological groups over $\S$}
$$0\to J_{\S}(\H)^0\to J_{\S}(\H)\to G\to 0.$$

\medskip\noindent
{\bf 2.5.~Remarks.}
(i) The topological group $G$ over $\S$ is locally homeomorphic to
$\S$ (on a neighborhood of each point of $G$).
However, its topology {\it cannot be} Hausdorff unless $|G_s|=1$
for all points $s$ of $\S\setminus S$.
Moreover, it is unclear whether the restriction $G|_{S_{\a}}$ of
$G$ over $S_{\a}$ is locally trivial over $S_{\a}$ (on a
neighborhood of each point $s$ of $S_{\a}$).
It is quite difficult to determine $G_s$ in general.

(ii) In case $\H$ is pure of weight $-1$, it is known
(see e.g.\ [1]) that there is an inclusion
$$G_s\otimes_{\Z}\Q\into\Hom_{\MHS}(\Q,H^1i_s^*j_{!*}\H_{\Q}),
\leqno(2.5.1)$$
where $j_{!*}\H$ denotes the intermediate direct image.
If $D=\S\setminus S$ is a divisor with normal crossings, then the
target is calculated by a subcomplex of the Koszul complex as is
well-known.
If furthermore $\H$ is a nilpotent orbit or if $\H$ corresponds to
a family of Abelian varieties, then it is easy to show the surjectivity of (2.5.1) although it does not hold in general.
(This subject will be treated in a forthcoming paper.)

\bigskip\bigskip
\centerline{{\bf References}}

\medskip
{\mfont
\item{[1]}
P.~Brosnan, H.~Fang, Z.~Nie and G.~Pearlstein,
Singularities of admissible normal functions (preprint
arXiv:0711.0964).

\item{[2]}
J.~Carlson, Extensions of mixed Hodge structures, in Journ\'ees de
G\'eom\'etrie Alg\'ebrique d'Angers 1979, Sijthoff-Noordhoff Alphen
a/d Rijn, 1980, pp.~107--128.

\item{[3]}
H.~Clemens, The Neron model for families of intermediate Jacobians
acquiring ``algebraic'' singularities, Publ.\ Math.\ IHES 58 (1983),
5--18.

\item{[4]}
P.~Deligne, Equations diff\'erentielles\`a points singuliers
r\'eguliers, Lect. Notes in Math. vol. 163, Springer, Berlin, 1970.

\item{[5]}
P.~Deligne, Th\'eorie de Hodge II, Publ. Math. IHES 40 (1971),
5--58.

\item{[6]}
F.~El Zein and S.~Zucker, Extendability of normal functions associated
to algebraic cycles, in Topics in transcendental algebraic geometry,
Ann. Math. Stud., 106, Princeton Univ. Press, Princeton, N.J., 1984,
pp.~269--288.

\item{[7]}
M.~Green, P.~Griffiths and M.~Kerr,
N\'eron models and limits of Abel-Jacobi mappings (preprint).

\item{[8]}
M.~Kashiwara, A study of variation of mixed Hodge structure,
Publ.\ RIMS, Kyoto Univ. 22 (1986), 991--1024.

\item{[9]}
M.~Saito, Mixed Hodge modules, Publ.\ RIMS, Kyoto Univ. 26
(1990), 221--333.

\item{[10]}
M.~Saito,
Admissible normal functions, J.\ Algebraic Geom.\ 5 (1996),
235--276.

\item{[11]}
W.~Schmid, Variation of Hodge structure: The singularities of the
period mapping, Inv.\ Math.\ 22 (1973), 211--319.

\item{[12]}
J.H.M.~Steenbrink and S.~Zucker, Variation of mixed Hodge structure,
I, Inv.\ Math.\ 80 (1985), 489--542.

\item{[13]}
S.~Zucker, Generalized intermediate Jacobians and the theorem on
normal functions, Inv.\ Math.\ 33 (1976),185--222.

\medskip
{\sfont
\baselineskip=10pt
Department of Mathematics, The University of British Columbia,
1984 Mathematics Road,

\quad
Vancouver, Canada

Department of Mathematics, Michigan State University,
East Lansing, MI 48824 USA

RIMS Kyoto University, Kyoto 606-8502 Japan

\smallskip
\vers
}}
\bye